\documentclass[12pt, a4paper,  reqno, english]{amsart}
\usepackage[all]{xy}
\usepackage[margin=0.9in]{geometry} 
\usepackage{amssymb,amsmath,amsthm}
\numberwithin{equation}{section}
\usepackage{lineno}

  \usepackage{pgfplots}
\pgfplotsset{compat=newest}
\usepackage{tikz} 
\usetikzlibrary{matrix,decorations.pathreplacing}
\usetikzlibrary{calc,angles,positioning,intersections,quotes,decorations.markings}
\usepackage{tkz-euclide}

\theoremstyle{definition}

\newcommand{\ben}{\begin{enumerate}}
\newcommand{\een}{\end{enumerate}}
\newcommand{\eit}{\begin{itemize}}
\newcommand{\beq}{\begin{equation}}
\newcommand{\eeq}{\end{equation}}

\renewcommand{\leq}{\leqslant}
\renewcommand{\geq}{\geqslant}

\begin{document}

\title{It is not ``B\'ezout's identity''} 
\author{Andrew Granville}

 \thanks{Thanks to Fran\c cois Lalonde, Karen Parshall and John Voigt for helpful conversations and emails. Special thanks to Ken Saito for insightful remarks and references about Euclid and Greek mathematics, Louise Poirier for clarifications about Inuit mathematics, and Emmanuel Kowalski and Jean-Pierre Serre for correspondence about the history of Bourbaki.}
 \dedicatory{A rose by any other name would smell as sweet --- \textsc{William Shakespeare}}
 \begin{abstract}  Given two non-zero integers $a$ and $b$ there exist integers $m$ and $n$ for which 
\[
am-bn = \text{gcd}(a,b) .
\]
An increasing number of  mathematicians have been  calling this   ``B\'ezout's identity'', ``B\'ezout's lemma'' or  even ``B\'ezout's theorem'', some encouraged 
 by finding ``identit\'e de B\'ezout''   in Bourbaki's \emph{\'El\'ements de math\'ematique}.
Moreover the observation that if $\gcd(a,b)=1$ then this is an ``if and only if'' condition, is sometimes called the ``Bachet-B\'ezout theorem''. 

However (we will explain that) all of this is in Euclid's work from around 300 B.C., when his writings are interpreted in context. So why does he not get credit?  Some authors learned the name ``B\'ezout's identity'' and have perhaps not consulted Euclid, so copied the misattribution.  Others, like some Nicolas Bourbaki collaborators, have perhaps browsed Euclid's results, but in a form written for the modern mathematician, and missed out on what he really did (though certainly others, such as Weil, did not). In this article we will carefully explain what Euclid's arguments are and what his approach was. We will also share 
 Emmanuel Kowalski's guess as to the reasons behind  Bourbaki's misnomer.

To appreciate Euclid, you need to read his work in appropriate context: Lengths were the central object of study to the geometer Euclid, though  he brilliantly developed the theory of the numbers that measured those lengths. Today's mathematicians read his number theory results as being about abstract  numbers not measurements. However the correct interpretation changes how these results are perceived; moreover Euclid's proofs make clear Euclid's intentions. 

These misperceptions reflect recent discussions about the difficulties sometimes faced by indigenous people when learning mathematics. We will discuss how some indigenous groups may learn numbers in certain practical contexts, not as abstract entities, and struggle when curricula assume that we all share abstract numbers as a basic, primary fully-absorbed working tool.
 \end{abstract}

\maketitle

\section{The genesis and use of Euclid's algorithm in Euclid's \emph{Elements}}

Most readers will have understood an abstract notion of a number as a small child, learning to count and then appreciating that numbers can be applied to all sorts of different situations. This seems so natural that one can scarcely imagine doing otherwise and yet early societies (like the ancient Babylonians or Greeks) and isolated societies today (like indigenous groups around the world) could proceed quite differently in educating their infants.   

We don't know how ancient Greeks learned about numbers as children but one would guess that a reader of the   \emph{Elements}  would have been encouraged, by their culture, to take a geometric perpective.
Indeed when studying ancient texts like Euclid's \emph{Elements}, as well as the existing fragments of its predecessors, it seems apparent that well-educated ancient Greeks working with proofs, were first steeped in geometric concepts and constructions, then found the need to measure and compare lengths, and so developed their understanding and use of (positive) integers. 

\subsection{Throughout history, men have compared lengths} 
In around 590 B.C., Thales of Miletus was said to have calculated the height of a pyramid by comparing the length of its shadow with the shadow of a small stick, a ``\emph{gnomon}'', noting that they subtend similar triangles and so that the two ratios of Height to Length are equal. Indeed  a gnomon dating from 
as far back as 2300 BC has been excavated in China, and they were probably then in common use.

Today we measure an angle  as a fractional multiple of $2\pi$ (or of $360^o$).  Ancient Greek sources would (equivalently) measure an angle by presenting it as the ratio of the  circumference of the circle to the arc on the circle subtended by the given angle. This had important applications: for example, ancients   realized that the angle that the shadow of a gnomon subtends at midday would tell one how far one is from the perpendicular to the sun.\footnote{If one is at the perpendicular then the gnomon will have no shadow as the sun shines directly down; a shadow emerges as one moves away from the perpendicular.}  Moreover one can use such measurements to estimate the \emph{axial tilt} which is the angle between the tangent plane to the line between the centres of the earth and the sun, and the axis of rotation of the earth.  In around 200 B.C., 
Eratosthenes  gave the approximation $83:11$ for the ratio  of a half circumference to   the axial tilt (that is, the angle is about $\frac{11\pi}{83}$ in modern parlance).\footnote{This was also estimated even earlier by Pytheas in the Greek colony in Massilia (modern day Marseilles) in around 350 B.C.}    Ptolemy (around 140 A.D.) suggests Eratosthenes obtained this ratio using the method of \emph{anthyphairesis}, a precursor to  the Euclidean algorithm,  which we now explain.

\subsection{Aristotle's \emph{anthyphairesis}} This is mentioned in the  ``\emph{Topics}'' of Aristotle (384-322 BC)\footnote{One of his six books on logic collectively known as ``\emph{The Organon}''.} and was used to determine when possible,  otherwise to approximate, the ratio between two given lengths. Suppose for now the lengths are given by two pieces of string, one longer than the other. The idea is to ``remove'' the length of the smaller from the larger, and then repeat, over and over. 
We start with two pieces of strings pulled tight, AB of length $x$, and CD of length $y$ with $x<y$.  
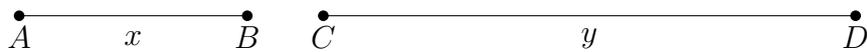
 \begin{center}
\begin{tikzpicture}
    \draw (0,0) -- (3,0);
    \node at (1.5,-0.3) {$x$};
    \fill (0,0) circle (2pt) node[below] {$A$};
    \fill (3,0) circle (2pt) node[below] {$B$};
    \draw (4,0) -- (11,0);
    \node at (7.5,-0.3) {$y$};
    \fill (4,0) circle (2pt) node[below] {$C$};
    \fill (11,0) circle (2pt) node[below] {$D$};
\end{tikzpicture}
\end{center}
Now we square up the far ends of the two pieces of string (B and D in the diagram below), and then the length from A to C is the difference $y-x$:
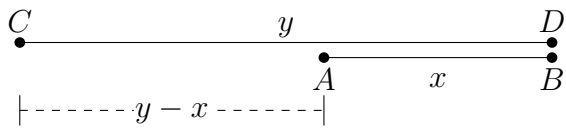
\begin{center}
\begin{tikzpicture}
    \draw (4.5,0) -- (7.5,0);
    \node at (6,-0.3) {$x$};
     \fill (4.5,0) circle (2pt) node[below] {$A$};
    
    \draw (.5,.2) -- (7.5,.2);
    \node at (4,0.4) {$y$};
      \node at (2.5,-0.7) {$y-x$};
        \fill (.5,0.2) circle (2pt) node[above] {$C$};  \fill (7.5,0.2) circle (2pt) node[above] {$D$};
    \fill (7.5,0) circle (2pt) node[below] {$B$};
    
           \draw[dashed] (.5,-.7) -- (2,-.7); \draw[dashed] (3.1,-.7) -- (4.5,-.7);
       \draw (.5,-.5) -- (.5,-.9);  \draw (4.5,-.5) -- (4.5,-.9);
           \end{tikzpicture}
\end{center}
The ancient Greeks did not necessarily label the lengths of the pieces of string (our $x$ and $y$), sometimes only their endpoints, as their actual lengths were not important, only the ratio of their lengths.

In our last diagram it is evident that we could have removed the length of AB from the distance from C to A a second time and still get left with something positive:
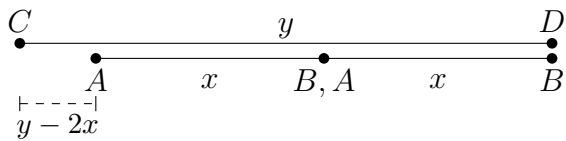
\begin{center}
\begin{tikzpicture}
    \draw (4.5,0) -- (7.5,0);
    \node at (6,-0.3) {$x$};
     \fill (4.5,0) circle (2pt) node[below] {$B,A$};
      \draw (1.5,0) -- (4.5,0);
    \node at (3,-0.3) {$x$};
     \fill (1.5,0) circle (2pt) node[below] {$A$};
    
    \draw (.5,.2) -- (7.5,.2);
    \node at (4,0.4) {$y$};
      \node at (1,-.9) {{\small $y-2x$}};
        \fill (.5,0.2) circle (2pt) node[above] {$C$};  \fill (7.5,0.2) circle (2pt) node[above] {$D$};
    \fill (7.5,0) circle (2pt) node[below] {$B$};
    
           \draw[dashed] (.5,-.6) -- (1.5,-.6);
       \draw (.5,-.5) -- (.5,-.7);  \draw (1.5,-.5) -- (1.5,-.7);
           \end{tikzpicture}
\end{center}
The length that is now left, $y-2x$, is shorter than $x$, the length of AB, so we cannot remove its length again.
Indeed no matter how much larger $y$ is than $x$, there is some integer number of times, ``the \emph{quotient}'', that we can remove $x$ from $y$ and no more.

We now repeat this algorithm working with the   length $x$ and the newly obtained smaller  length $y-2x$, again subtracting the smaller length from the larger, an appropriate number of times.
How is this useful?

Suppose that $x=15$ and $y=54$. We first subtract 3 times $15$ from $54$ to get $9$,\footnote{Here $3$ is the ``quotient'' and can be obtained as $\lfloor \frac {54}{15}\rfloor$.} so we next compare $9$ and $15$.
We subtract 9 once from $15$ to obtain $6$. We then subtract $6$ once from $9$ to obtain $3$, and $3$ twice from $6$ to obtain $0$.
Therefore the sequence of quotients we subtracted is $3,1,1,2$ to get to $0$.  If we go through the same process with $10$ and $36$ then we find the \emph{same sequence} $3,1,1,2$ of quotients. Aristotle would have then deduced that the ratios $15:54$ and $10:36$ are equal. We would get the same ratio if the lengths were $5e$ and $18e$, or indeed with any pair of starting numbers as long as the sequence of quotients is exactly the same, so this process does not need  whole numbers to work.  This process might seem cumbersome but \emph{this is} how mathematical thinking developed, and these quotients yield the \emph{continued fraction} for our ratio. That is,
\[
\frac{54}{15}= 3 + \frac 1{1 + \frac 1{1 + \frac 12}} ,
\]
though it is more convenient to write this as $\frac{54}{15}=[3,1,1,2]$

Given two arbitrary pieces of string it is unclear whether this anthyphairesis process will ever end, so one can get approximations by stopping the process after a certain number of steps. For example, if we attempt this with the ratio of the circumference of a circle to its diameter, then the quotients   are 
$3,7,15,1,292,\ldots$ and it is difficult to guess whether the process ends.\footnote{Especially as the lengths we are now working with are a tiny fraction of the original lengths, so difficult to deal with in practice.}  However each successive step of the algorithm gives rise to better and better approximations to our ratio (which we denote by $\pi$):
\[
3,\  3 + \frac 17=\frac{22}7,\ 3 + \frac 1{7 + \frac 1{15}}=\frac{333}{106},\ 3 + \frac 1{7 + \frac 1{15+\frac 11}}=\frac{355}{113}
\]
and this last approximation differs from $\pi$ by less than $10^{-6}$.

It was realized by the ancient Greeks that if this process comes to an end then the two lengths are \emph{commensurate}, that is they are a rational multiple of each other, and if not they are \emph{incommensurate}, which means their ratio is irrational.  Euclid's first main contribution to number theory was to develop the theory of commensurate lengths.

\subsection{Measuring lengths with two given measures}
Following in the Greek tradition, Euclid was interested in what lengths could be measured by two given lengths (we shall again assume these are lengths of string, but they could be measuring sticks).  The basic operation of removing the length of the smaller from the larger was given above; it is even easier to add the lengths by putting these strings end-to-end:

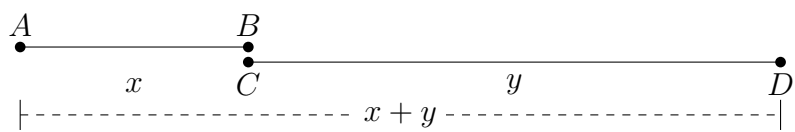
\begin{center}
\begin{tikzpicture}
    \draw (.5,0.2) -- (3.5,0.2);
    \node at (2,-0.3) {$x$};
     \fill (3.5,0) circle (2pt) node[below] {$C$};
      \fill (.5,0.2) circle (2pt) node[above] {$A$};
    
    \draw (3.5,0) -- (10.5,0);
  \node at (7,-.3) {$y$};
        \fill (3.5,0.2) circle (2pt) node[above] {$B$};
    \fill (10.5,0) circle (2pt) node[below] {$D$};
    
       \draw[dashed] (.5,-.7) -- (4.9,-.7); \draw[dashed] (6.1,-.7) -- (10.5,-.7);
       \draw (.5,-.5) -- (.5,-.9);  \draw (10.5,-.5) -- (10.5,-.9);
           \node at (5.5,-0.7) {$x+y$};
\end{tikzpicture}
\end{center}

\noindent Here we begin with AB, and then use the measure of length $y$ to go on from $B$ to $D$ so that $AD$  now has length $x+y$.

Euclid begins with a unit of length, and assumes that this can be used to measure each length of string, the length AB is $x$ times the unit length where $x$ is an integer (that is, it can found precisely by taking the unit length, end-to-end, $x$ times), and the length CD is $y$ times the unit length where $y$ is an integer.  

Euclid's algorithm is essentially the study of finite anthyphairesis processes. 
All lengths that occur through the anthyphairesis process are to be measured using the original strings of lengths $x$ and $y$, putting them end-to-end in one direction or the other. This was obvious to Euclid as his discussion is framed in these terms, but the statements of his results do not appear to make this point when read out of this context. Let's discuss this further with the example of lengths  $12$ and $41$. First we subtract $12$ three times from $41$ to obtain the  length $41-3\cdot 12=5$. Then we subtract   5 twice from $12$ to obtain the length $12-2\cdot 5=2$. To see that the  length $2$ can be measured by our original strings of lengths $12$ and $41$ we unpack this explicitly to get 
\begin{align*}
2&=12-5- 5=12-(41-12-12-12)-(41-12-12-12)\\
&=12-41+12+12+12-41+12+12+12,
\end{align*}
which is $7\cdot 12-2\cdot 41$.
We can go one step further subtracting 2 twice from 5 to get
$5-2\cdot 2=1$ so that
\[
1=(41-3\cdot 12)-2\cdot (7\cdot 12-2\cdot 41)=41-3\cdot 12-14\cdot 12+4\cdot 41=5\cdot 41-17\cdot 12.
\]
Complicated but doable. Euclid did not explicitly do this ``unpacking'', hampered by his notation and by giving no examples, but there can little question that this is what he knew himself to be doing.

In more modern language $5, 2$ and $1$ are each   $\mathbb Z$-linear combinations of $12$ and $41$. 

\subsection{Euclid's conclusions} 
Euclid established that the anthyphairesis algorithm always terminates with a length which could measure both of the original lengths. So in the example starting with $15$ and $54$ we terminated with $3$, and each of $15$ and $54$ are multiples of $3$. He also proved that any common measure of the two starting lengths would be a measure of the length obtained at the end of the algorithm, so the terminal length is actually the   \emph{greatest common measure},
the longest length that measures both $15$ and $54$.
(Today we would call this the   \emph{greatest common divisor}.)
Therefore his proofs yield the following:
\smallskip

\begin{quote}
\emph{Suppose we are  given positive integer lengths $x$ and $y$, with greatest common measure $z$. We can determine lengths $X$ and $Y$ by the anthyphairesis process, which differ by   $z$, where $X$ can be measured by $x$ (that is $X$ is a multiple of $x$), and  $Y$ can be measured by $y$.}
\end{quote}
\smallskip

\noindent In modern conventions, he proved that there exist non-negative integers $m$ and $n$ for which
\[
|mx-ny| = \text{gcd} (x,y)
\]
(where $X=mx, Y=ny$ and $z=\text{gcd} (x,y)$).

The terminal length, the greatest common measure, obviously measures any multiple of its length, so Euclid also knew the following:
\smallskip

\begin{quote}
\emph{Suppose we are  given positive integer lengths $x$ and $y$, with greatest common measure $z$. For any length $Z$ that can be measured by $z$,  we can determine lengths $X$ and $Y$ by the anthyphairesis process, which differ by   $Z$, where $X$ can be measured by $x$, and  $Y$ can be measured by $y$.}
\end{quote}
\smallskip

\noindent In modern mathspeak, we drop the whole notion of measures and lengths and deal solely with the integers that represent those lengths.
Thus Euclid's results can be formulated as follows:

\begin{quote}
\emph{ For two positive integers $x$ and $y$  with greatest common divisor $z$ we have
\[
\{  |ax-by| :  a,b \in \mathbb Z_{\geq 0}\} = \{  cz :  c \in \mathbb Z_{\geq 0}\}.
\]
In particular, there are non-negative integers $m$ and $n$, which may be found using the Euclidean algorithm, for which 
$mx-ny=z$ or $ny-mx=z$.}
\end{quote}
\smallskip

\noindent To go from here to obtain the standard results of the modern number theory curriculum we have to  include negative integers in our formulation and conclusion. These did not arise in Euclid's formulation since he was dealing with length and for the ancient Greeks there was some discomfort in working with negative quantities since they did not as evidently represent a physical entity as positive quantities.  Nonetheless, today this is a breeze:

\begin{quote}
\emph{ For two non-zero integers $x$ and $y$  with greatest common divisor $z$ we have
\[
\{ ax+by :  a,b \in \mathbb Z\} = \{  cz :  c \in \mathbb Z\}.
\]
In particular, we can determine integers $m$ and $n$ using the Euclidean algorithm, for which }
\[
mx-ny=z.
\]
\end{quote}
\smallskip

I hope the above makes clear Euclid's intent in his writing. Although he replaced lengths by positive integers in the formulation of his results (as discussed below), Euclid did describe divisibility using the concept of measuring, which incorporates more meaning in the    context of his time and his thinking.
Any interpretation of ancient literature should surely take into account the extant literature, protocols  and notation at the time of the writings.
What I have proposed is what I believe Euclid meant, though my  explanation is consciously tainted by an anachronistic understanding of where this all leads, for example by my use of ``$\mathbb Z$-linear combinations'' and by my desire to end up phrasing Euclid's results in terms of the ideal generated by $x$ and $y$ over $\mathbb Z$.

\section{Interpreting Euclid's text}

Let me now tie in the above discussion with the appropriate statements and proofs in Euclid's actual writings:
The paragraphs corresponding to Propositions 1 and 2 and the subsequent ``Porism'' in Book VII of Euclid's \emph{Elements} make this argument as well as several deductions.  Here we quote some of the relevant parts of Heath's  translation:

\subsection*{Definition} 
\emph{The greater number is a \emph{multiple} of the less when it is measured by the less.}\footnote{Euclid gave 22 definitions at the start of this book. Indeed, in any math book the author has to judge   what level of mathematical knowledge  one can assume? What might be too self-evident to define?   Euclid does not define ``measure'', and he does not really define addition, presumably because he felt these were common knowledge, but these are glaring omissions from a modern axiomatic perspective.}\smallskip

Here, Euclid defines a ``multiple'' as one number, which gives the length of a larger stick, being measured exactly by a smaller stick.
This is not the same as the modern more abstract definition ($b$ is a multiple of $a$ if there exists an integer $k$ for which $b=ak$), in particular because Euclid never mentions the quotient $k$.  Euclid frequently discusses ratios of lengths (see section \ref{sec: ratios} below), so  in his mathematics   the ratio $b:a$ is the same as the ratio $k:1$ for some integer $k$, but this evidently is a complicated formulation, and so largely avoided unless necessary.
 
 To Euclid and his intended readers it is natural to measure lengths, and ratios of lengths. Introducing the abstract notation of a non-negative integer in the  statement of results reads like it was slightly unusual, something that needed thinking about. To the modern era mathematician the opposite is true; we are taught the abstract notion of a non-negative integer (and what it might represent) at an early age, and it is a simple application to think of an integer as measuring a length. For us, it is unusual to use measuring lengths as an agreed upon base for proving number theory theorems.

\subsection*{Proposition 1} \emph{Two unequal numbers being set out, and the less being continually subtracted in turn from the greater, if the number which is left never measures the one before it until an unit is left, the original numbers will be prime to one another.}\smallskip

In this statement of Proposition 1, Euclid describes somewhat tersely the anthyphairesis process for successively subtracting the smaller integer from the larger, what we now call  the ``\emph{Euclidean algorithm}''  when applied to a pair of positive integers.\footnote{It is valuable to note that the anthyphairesis algorithm as described above only assumed that the two starting numbers are positive real numbers.}  He then observes that if at the end one gets $1$ then the original two numbers have no common factor.  

Euclid appears to have written this statement for abstract integers without associating them to  lengths, especially if one does not take care with the intended meaning of the word ``measure'', and that seems to be how it is usually (mis-)interpreted.
Yet in his proof it is evident that he did mean lengths. Indeed his proof begins:\footnote{This separation into statement and proof is a more modern concept. In Euclid's text, a proof is not given a name, separating it from the statement of the result. It is simply the next paragraph.} 
\begin{quote}
\emph{For, the less of two unequal numbers AB, CD being continually subtracted from the greater, let the number which is left never measure the one before it until an unit is left; I say that AB, CD are prime to one another, that is, that an unit alone measures AB, CD.} 
\end{quote}
So we immediately see that a number is represented as the distance between two points (AB, and then CD). 

Above we saw that in Euclid's formulation that any length represented by the anthyphairesis process is evidently an integer linear combination of the original two lengths ($x$ and $y$). Therefore when Euclid writes that the process ends in a unit, what he is thinking is that there are lengths $X$, measured by $x$, and $Y$ measured by $y$, for which $|X-Y|=1$.  He then deduces that if this the case then $\gcd(x,y)=1$. This is therefore the converse part of the so-called ``Bachet-B\'ezout theorem''.  However to clarify this it would be helpful to prove that the Euclidean algorithm ends with $\gcd(x,y)$. This is what Euclid proves next:

\subsection*{Proposition 2} \emph{Given two numbers not prime to one another, to find their greatest common measure.} \smallskip

In this result Euclid uses the anthyphairesis process to find the gcd of two integers (whether or not they are prime to one another).
As an abstract statement, interpreted out of this context, this does appear to only say that the Euclidean algorithm ends with $\gcd(a,b)$.
However the proof \emph{immediately} yields the more useful formulation in terms of multiples. We now give part of Euclid's proof in a little more detail:

In the proof of Proposition 2 Euclid supposes that CD does not ``measure''  AB (that is, a stick from C to D, a stick of length $|CD|$, cannot be used to 
``measure'' from A to B, a stick of length $|AB|$). He then lets EB be the longest part of AB that  can be measured by CD (i.e that is, EB represents the largest integer  multiple of the length of CD that is no longer than AB), and so we get a left-over piece AE, which is less then AB in length. Euclid notes that this is ``measured'' by AB and CD (that is, $|AE|$  is a $\mathbb Z$-linear combination of $|AB|$ and $|CD|$). He then repeats this argument: If AE doesn't measure CD then let FD be the longest part that is measured by AE, and then this leaves CF.  Again we see that CF is given by the length CD and then removing the correct number of multiples of AE, so $|CF|$ is a $\mathbb Z$-linear combination  of $|CD|$ and $|AE|$.  

 Euclid does not then explicitly state that, since $|AE|$ is a $\mathbb Z$- linear combination of $|AB|$ and $|CD|$, we can deduce that $|CF|$ is a $\mathbb Z$-linear combination  of $|AB|$ and $|CD|$ (as in my strings of lengths 12 and 41 example above), though I believe he would have felt this was obvious in context; that is, numbers being used to measure, not as abstract entities. He is more explicit about this in his proof of Proposition  1 though his goal there is slightly different, in that he is looking at the consequences when $a$ and $b$ have no common factor other than 1. However I think that his intent throughout the proof of Proposition  2 is clear; he is looking at what lengths can be measured by the original two strings.

Euclid next supposes that CF measures AE (that is, the length of AE is a multiple of the length of CF). Since AE measures FD, he deduces that CF measures FD, and so measures CF$\cup$FD$=$CD.  Now CD measures EB and so CF also measures EB.  But CF also measures AE so it  measures the whole of AB.

  Phew, he is done, he has proved that the last obtained measure, the length of CF, indeed measures the two original  lengths, of AB and of CD. In other words, he proves that the (single) terminal integer length which is the  outcome of the algorithm divides the two original integer lengths.
  
   Next he goes on to prove that this terminal  length  is indeed the \emph{largest} common measure (that is the algorithm terminates with gcd$(x,y)$), and in the ``Porism'' he notes that the same proof implies that any common measure divides the largest common measure (that is, all common divisors of $x$ and $y$ divides gcd$(x,y)$):

\subsection*{Porism (Corollary)} \emph{ If a number measure two numbers, it will also measure their greatest common measure.}\smallskip

 This reaffirms Proposition 1 in the special case $\gcd(a,b)=1$. So we see that Euclid, correctly interpreted, indeed proved all of the results that some of the literature claims for much later French authors!
\smallskip

\section{A different way to present ideas} 
The modern mathematician might wonder why Euclid did not run his algorithm to an arbitrary number of steps (as any textbook today would), rather than stop at the second step.\footnote{Indeed this does not complete the proof for our lengths 12 and 41 example above.} At that time there was no subscript notation,   no ``$\ldots$'', nor even $\& c$ (as Gauss used in 1804) nor ``et cetera'' (which is a 14th century invention).  Up to that time the reader would have been expected to understand that the exactly analogous argument would work if one of the two lengths divided the other at the third step, or the fourth, or the fifty-seventh, just as we do today, though we now have several notations to describe that. Thus Euclid was not ``wrong'' but rather writing for the (educated) reader of his time who would have had a feeling for what was meant in this case and others. Similarly I believe that a contemporary reader of Euclid felt, like him, that he was writing about measuring lengths, and would not have forgot this halfway through the proof, so even if he failed to explicitly mention it,   his audience would have quickly caught on.
 
 But someone might complain that here I have been talking about proofs and not statements of theorems so surely that is what we should use in interpreting what Euclid meant? However, Euclid does not delineate results in this (modern) way. He does usually write a starting sentence and ending sentence in the section dedicated to each proposition, describing what he is about to do, and then what he has done (and how it can be interpreted) but, again, I think it is clear that for him a number measured something, so his more abstract first statements (when interpreted in a modern way) should surely be interpreted in the sense that he meant them.\footnote{Euclid never sullied Book 7 with an example or a specific integer (other than 1). My feeling is that this is because he was interested in ratios (what lengths measured what other lengths) so he consciously wanted to deter his readers from choosing non-generic units.}

We cannot easily get inside Euclid's head because what he might have seen as obvious from the information he had at hand, does not necessarily fit the way we construct mathematics in the curriculum today.  However it seems quite clear that the positive integers that he constructs in his algorithm are positive integers that he knows can be written as a linear combination of the original two integers; in modern parlance,\footnote{And not worrying about the anachronistic use of negative integers rather than using positive integers and absolute values.} for given  integers $a$ and $b$, he finds successively smaller positive integers in the ideal
\[
\{ am+bn: m,n\in \mathbb Z\}
\]
and proves that the smallest positive integer in the ideal is the \emph{greatest common divisor} of $a$ and $b$.  The relatively recent formalism of ideals has been very useful.
I am not suggesting for one moment that Euclid would have guessed at this but subsequently one sees that his work leads naturally to this concept that revolutionized 19th century algebra.\footnote{To be more precise, the revolution came from discovering that not all ideals are principal (that is,  in number fields they are not necessarily   generated by one element like a gcd).}   The relevant developers of these ideas, Dirichlet, Kummer and Dedekind, would have known their Euclid (and not read it in translation!), and would not have misinterpreted his meaning.

\section{Diophantus} Only two ancient Greek books have had a direct impact on modern number theory, both authors from Alexandria, a Hellenized city of ancient Egypt. The first was Euclid's \emph{Elements} from around 300 B.C., the second Diophantus's ``\emph{Arithmetica}'' from around 260 A.D.
While Euclid's work was an attempt to axiomatize and formally prove the basics of geometry and consequently number theory, Diophantus was interested in ad hoc problems, finding positive integer and rational solutions to various explicit equations, some of higher degree. One might have hoped that these solutions lead  to some general principles but this is not obvious, certainly as written.

There were a lot of changes in Greek society  in the intervening 500 years between Euclid and Diophantus, most of which we can only guess at.   We do know that Euclid's scholarship was very much appreciated in his time, whereas Diophantus's much less so.  

For this article we are primarily interested in how Diophantus approached and interpreted integers:  It is clear that he treated them as abstract entities in his writings referring to an unknown integer as 
``a number of units which is undefined''. He takes powers and in doing so uses abbreviations (though his notation is not very logical compared to today's notation for powers). 
Taking powers was really not possible for Euclid since multiplication was not defined in a convenient way (if a length $a$ measured $b$ then Euclid might make explicit how many copies of $a$ go into making $b$, but the roles of $a$ and this quotient were only defined in those terms, and they certainly were not treated as interchangeable).

Diophantus's notion of integers still had some restrictions that seem left over from an earlier era, in that he did not allow negative integers or even zero. On the other hand he would subtract one quantity from another but it was unacceptable if the consequence was negative. Here are some details (simplified from the examples in his text):

--- He would  treat finding integer solutions to $x^2+x=y^2$ and finding integer solutions to $x^2-x=y^2$ as distinct problems. 

--- When solving quadratic equations, he would only take the positive square root. For example when solving $x^2+3=4x$ using the usual quadratic formula, he would have only acknowledged the solution $x=3$, not $x=1$ since, in modern notation, the two  solutions are $2\pm \sqrt{1}$.

--- When solving $x^2=x$ he would immediately divide through by $x$ to obtain $x=1$, ignoring the solution $x=0$. 

We do see a progression of thought from treating integers as a way to quantify length to becoming an abstract entity, albeit not developed as we would today.
I do not know how we got to the modern abstract notion of integers, and when this seemed like the natural way to proceed.

\section{Re-discovery}
Before fairly modern times, books and papers were not so easily accessible, and indeed  mathematicians who were just starting out might not have known where to look for previous work on ideas that they were developing.  And so they might re-discover the works of the masters, which can be off-putting to an ambitious novice, but they might also have brought a new, unexpected approach to the subject.

\subsection{Bachet}  Bachet de M\'eziriac (1581-1638)  is an immensely important mathematician because his translation of, and commentaries on, Diophantus's ``\emph{Arithmetica}'' from ancient Greek  into Latin in 1621, made it broadly accessible in the Renaissance, and inspired a new generation of mathematicians. For example, a young Fermat defaced his copy of Bachet's translation with his infamous marginal note, a remark that became known as ``Fermat's Last Theorem''.
Diophantus's book  re-inforced the study of integer and rational solutions to equations as central to our subject.

Bachet had already published a book \cite{Ba} in 1612 of  primarily arithmetical, mathematical puzzles.\footnote{For example, ``questions involving number bases other than 10; card tricks; watch-dial puzzles depending on numbering schemes; think-of-a-number problems; river crossing or ferry problems ...; problems concerning magic squares...; the Josephus problem...; various weighing problems; and liquid pouring problems'' (taken from the online MacTutor History of Mathematics Archive).}
In the second edition, published in 1624, Bachet determined the general solution to $ax-by=\pm 1$ (\cite{Ba2}, Proposition XVIII).  It is evident that Bachet had not read Euclid's work as he laboriously discovers for himself a crude version of the Euclidean algorithm. However, unlike Euclid, Bachet knows to treat integers directly (with no mention of lengths or measures), he works with variables to represent the integers without the need for explanation and  he   multiplies and divides assuming his readers can appreciate what he is doing without further explanation. By today's standard his notation is limited. It is noticeable that he includes a non-trivial numerical example in a box on every page or two to help the reader's understanding.\footnote{And he does use far more variables than are really necessary for his proof.}

What impressed the mathematicians of the Renaisance (some of whom had read Euclid), was that 
Bachet not only found a first solution to $ax-by=\pm 1$, like Euclid, but also determined all solutions. Technically this  is a simple development if one already appreciates Euclid's work, but conceptually it is a massive leap, as it fully answers a general question of Diophantus type.
 However there seems little doubt that this was already resolved
by  the  Indian scholar $\overline{\text{A}}$ryabhata in around 500 A.D. (and feasibly  by Archimedes in Syracuse (Greece) in 250 B.C., possibly by certain ancient Chinese scholars) but certainly Bachet's solution was the first known to Renaissance mathematicians, the main source of our intellectual inheritance today.

Bachet  mentions his result in his (earlier) Diophantus translation. Of that,
Weil (chapter 1.V in \cite{We}) writes ``in 1621, Bachet, blissfully unaware (of course) of his Indian predecessors, but also of the connection with the seventh book of Euclid, claimed the same method emphatically as his own''.

\subsection{Gauss} In the introduction of his ``\emph{Disquisitiones Arithmeticae}'' (1804), Gauss wrote  that by the time he was 18 years old, ``the greater part of the first four sections had already been completed before I had seen anything of the work of the other geometers on this subject''.
 This includes his work on solving linear equations. In fact his approach is rather different from that of Euclid: Gauss develops the basic theory of congruences and so, for given coprime integers $a$ and $b$ he solves $am+bn=1$ (in section 27 of \cite{Disq}) by first showing that there is an integer $m$ in the range $1\leq m\leq b$  for which $am\equiv 1 \pmod b$ and then letting $n:=\frac{1-am}b$.\footnote{He had shown that such an $m$ exists by proving that $\{ ar: 1\leq r\leq b\}$ gives a complete set of residues mod $b$ and, in particular, there must be a solution   $m, 1\leq m\leq b$ to $am\equiv 1 \pmod b$.} He goes from there to find all solutions by noting that if
 $at+bu=1$ is any solution then $t\equiv m \pmod n$ and $u:=\frac{1-at}b$. Gauss \emph{then} goes on to explain how to find the first solution $m,n$ much more efficiently by using Euclid's algorithm.\footnote{In \cite{We}, Chapter 1, $\S IV$, Weil wrote ``The general method of solution for [$ax-by=m$] is essentially identical with the Euclidean algorithm'' which is patently untrue for Gauss's first proof, though is true for the proof of Bachet, since Weil is criticizing both Bachet and those who give credit to him and B\'ezout.}
 
 By the time Gauss was preparing his book for publication a few years later, he had discovered the extant literature, praising Euclid's Book 7, Diophantus, Fermat, Euler, Lagrange, Legendre, ``and a small number of others''  in his introduction  as ``having opened the door to this divine science, and discovered the inexhaustible wealth it contains''.  It is worth noting that Gauss approved of Euclid's ``elegance and rigour'' but (slightly) disparages Diophantus as his arguments, although requiring ``genius and penetration'',  ``are too particular and rarely lead to general conclusions''.
 
 In section 28 of \cite{Disq} Gauss states that ``Euler is the first that has given a resolution of these equations'' (that is, finding \emph{all integer solutions} to $am+bn=1$) by a method of ``substituting other variables for $m$ and $n$''; and noted that Lagrange in 1767 resolved this through continued fractions by writing
 $\frac ab=[c_0,c_1,\cdots,c_\ell]$ and then letting $\frac rs=[c_0,c_1,\cdots,c_{\ell-1}]$ to obtain $ar-bs =\pm 1$.\footnote{Although the ancient Greeks essentially constructed continued fractions for $\frac ab$, I do not know of a reference where they observed that the solution to $ar-bs=\pm 1$ could be found by truncating the continued fraction. This exhibits a benefit  to   formulating these ideas with numbers rather than lengths. To prove this works, first establish that if $\frac{p_n}{q_n}:=[c_0,c_1,\cdots,c_n]$   then  $\begin{pmatrix}  c_0&1\\1&0\end{pmatrix} \cdots \begin{pmatrix}  c_n&1\\1&0\end{pmatrix} =\begin{pmatrix}  p_n&p_{n-1} \\ q_n&q_{n-1}\end{pmatrix}$ (perhaps by induction). Taking determinants with $n=\ell$ we obtain $as-br=(-1)^{\ell+1}$.}
 
Books took a long time to typeset at the beginning of the 19th century, particularly mathematics books, so long after the early pages of \cite{Disq} were typeset (and perhaps printed), Gauss was able to revisit them. Thus at the end of \cite{Disq} we find a few revisions of the earlier chapters, entitled ``Additions of the author''. Here  Gauss writes that after reading  Lagrange's   ``Additions \`a l'Alg\`ebre d'Euler'', page 525, he discovered that Bachet (from the 16th century) deserves credit for first finding all solutions to $ax+by =\pm 1$, not  Euler.    

I find it surprising that Gauss gave credit to Euler then Bachet but not to Euclid. Gauss makes clear he read Euclid's Book 7 and found it elegant, so he must be distinguishing between finding one solution and finding them all.  Although finding all solutions is an interesting augmentation of Euclid's original theory, it is hardly challenging when phrased as a Diophantine problem, and certainly  compared to the much more subtle task of finding a first solution.
However Gauss's proof that there is a first solution is so simple (although highly ingenious)    that he might  have considered finding all solutions to be of equal depth and thus felt Bachet was more worthy than it would seem to those of us who began our understanding of number theory by appreciating the Euclidean algorithm.

\subsection{B\'ezout} 
B\'ezout (1730-1783) published ``\emph{Th\'eorie g\'en\'erale des \'equations alg\'ebriques}'' in 1779, in which he developed the theory of resultants and, most famously, found that there are $dD$ intersections of generic planar curves of degrees $d$ and $D$, a key starting point for the development of algebraic geometry.

In the same book, B\'ezout established that there is an analogy of  Euclid's algorithm for polynomials, reducing the degree of the larger degree polynomials at each step, in place of reducing the   size of the larger integer in Euclid's original algorithm. Not only is this an important observation in of itself but it also hints at a much bigger picture, that such a result should hold in some generality for example for all Euclidean domains (an integral domain for which the Euclidean algorithm works).

\subsection{Bourbaki} The analogous result is formulated and proved   for principal ideal domains (which include all Euclidean domains) in  Bourbaki's \cite{Bo} \emph{Alg\`ebre}, chapitre 7. It is here that Bourbaki bizarrely writes\footnote{A bizarre attribution and bizarre that ``B\'ezout'' is misspelled as  ``Bezout''.}

\centerline {Th\'eor\`eme 1 (``identit\'e de Bezout'').}
\smallskip

However there is no part of this result that is due to B\'ezout (apart from the inspiration to generalize Euclid)!
I don't know whether this attribution  appeared earlier in the literature but it certainly appears in abundance today, including Wikipedia pages in both English and French (the latter citing Bourbaki as the authority), and  in many textbooks, such as \cite{Has, Haw, Hill, JJ, Sam}.
Some authors try to improve the historical accuracy, for example
Colmez \cite{Co} (footnote on p12) ``this is due to Bachet (1624); B\'ezout .. showed the analogous result in the ring $K[X]$''
and Tenenbaum  \cite{Ten} (ex 12 on p23)   ``Bachet's Theorem (1624), better known, wrongly, as B\'ezout's theorem'' but as we have seen, Bachet also played no part in first proving this ``theorem'' or  ``identity''.

\subsection{Kowalski in the Bourbaki archives}

As Kowalski pointed out to me,\footnote{By email, when provoked to defend the honour of Bourbaki.} Bourbaki was well aware of who did what, and lauds Euclid's
arguments very highly   (``on ne peut qu'admirer la finesse et la
s\^uret\'e logique qui s'y manifestent"), mentioning that his arguments
remain pretty much those used in modern times (``ne diff\`erent gu\`ere en
substance de ceux du chap. VI, sec 1''). Moreover  Bourbaki also mentions that Hindu mathematicians before the 13th
century had methods to deal with arbitrary systems of integral linear equations.  So it seems unlikely that  the Bourbaki collaborators really wanted to attribute results of this kind to either Bachet or B\'ezout.  Indeed Weil was an integral part of Bourbaki when this book was written (late 1940s/early 1950s) and in his History book 
\cite{We} writes that Fermat and Wallis certainly read Bachet but ``surely knew their Euclid too well not to recognize the Euclidean algorithm there''.

Kowalski (11/05/2024) kindly scoured the Bourbaki archives\footnote{http://sites.mathdoc.fr/archives-bourbaki/PDF} to try to find out how this attribution got in to the eventual manuscript:  

--- On page 62 of the first draft (076$\_$iecnr$\_$084.pdf),  Dieudonn\'e writes (in translation): ``If $a$ and $b$ are two relatively prime integers, there exist two integers $p$ and $q$ such that 
 
\noindent (5) \centerline{ $1=pa+qb$}

\noindent (Bezout's identities).'' 

--- Corollary 1 on page 58 of the next (undated and unattributed) draft (077$\_$iecnr$\_$085.pdf) gives  the same statement though  $a$ and $b$ are now in a PID. However B\'ezout is not mentioned.

--- The Corollary   on page 28 of the next version (121$\_$nbr$\_$028.pdf) was written by Weil in December 1949 and more-or-less gives the result in the final published form: ``A finite set of elements $x_1,\dots,x_m$ of a PID   have no common factor if and only if there exist integers $a_1,\dots ,a_m$ for which $a_1x_1+\cdots + a_mx_m=1$.'' Again B\'ezout is not mentioned.

--- On page 11 of the final draft (140$\_$nbr$\_$043.pdf) the statement is only slightly modified, though it is now upgraded to being ``Theorem 1''. The attribution to B\'ezout re-appears, but now as ``(``identit\'e de Bezout'')''. It is unclear who wrote this final draft though my guess is Dieudonn\'e because of the same spelling mistake in B\'ezout's name.\footnote{Kowalski points out that both Bezout and 
Meziriac, without the ``\'e''s, appear on the front page of different early editions of their books, so perhaps not too much should be read into these spellings.}

\subsection{Serre's recollections} Serre was a member of the Bourbaki collective by the time of Weil's draft so I asked him (by email) whether he recalled what had happened. He replied on 13/05/2024 (translated):
``I don't remember any discussions in Bourbaki on this subject. Generally speaking, issues of a historical nature were not discussed in the meetings (fortunately, as this would have wasted a lot of time); they
were left to the specialists: Weil first, and Dieudonn\'e, then also
Delsarte and Chevalley. Much later, other members of Bourbaki, like Samuel, Cartier, Dixmier, Borel and I became interested in it.''

Serre summed up his feelings on this issue as follows: ``You are right to say that the result was known to Euclid, other than for notation, in the case of integers. It became known for polynomials, thanks to B\'ezout. Lots of terminologies are less justified than that.''

Reading this Kowalski (14/05/2024) wrote to me: ``I also agree with his [Serre's] point that
much of mathematical terminology is even less reasonable than B\'ezout's
identity, although Bourbaki always strained to be absolutely impeccable
on that point, usually by not giving attributions."

\subsection{Why did Bourbaki persist in this mis-attribution?}

Kowalski (05/05/2024) pointed out that very few
theorems are given a specific attribution in Bourbaki,
``so the quotes indicate to me that this is just a convenient designation, more
convenient than ``Euclid's identity" maybe or ``Euclid's theorem".  It is
indeed used a few times later on.''

In the same email Kowalski doubts that this attribution went to press without Weil's explicit approval 
and that Weil was perfectly aware that these results  are all due to Euclid in substance (as discussed above).  
So ``it might even be that this is an inside joke on Weil's part, aware that
some people attributed this to B\'ezout.'' In other words, this was Weil being sardonic, perhaps frustrated with how B\'ezout was getting credit from lesser mathematicians who had not consulted the original sources, and perhaps this is why ``identit\'e de Bezout'' appears in quotation marks?

Looking through \emph{Alg\`ebre} chapters 4-7 I found that  Proposition 6 in Book 4.2 page 15 is ``formule d'interpolation de Lagrange'' without quotation marks in the text, Theorem 2 of Book 5, page 80 is attributed to Gauss without quotation marks, Theorem 3 of Book 6, page 250   to Euler-Lagrange without quotation marks. However in Proposition 6 of Book 4.1, page 8, we have 
``identit\'e d'Euler'' in quotation marks though it is referred to as ``la relation d'Euler'' on page 66 without quotation marks.\footnote{Surprisingly this Euler's identity is not $e^{i\pi}+1=0$, but rather that if $f(x_1,\dots,x_m)$ is homogenous of degree $m$ then $mf=\sum_i \frac{\delta f}{\delta x_i}$. The problem with convenient attributions to great mathematicians is that they have too many great theorems!} From this evidence, it is feasible that   quotation marks are some sort of Bourbaki protocol whose meaning has been forgotten, or perhaps there is no real meaning.

If it is true that Dieudonn\'e re-instated the attribution ``B\'ezout's identity'' as I suspect, then we must ask why? Did he think it got dropped in the intervening drafts by accident? Or did he grow up with this name and believe it belonged? One might have guessed that Dieudonn\'e would have wanted to do better.

Perhaps we will never find out the truth about this inappropriate attribution, and why it has spread so far and wide, but I would like to see it stamped out.

\section{Incommensurability and anthyphairesis} \label{sec: ratios}
 In Plato's dialogue, \emph{Theaetetus} claims that, developing ideas of his teacher Theodorus, he established that if the length of the side of a square is not an integer, but the area of the square is an integer, then the length of its side cannot be commensurate with 1, ``and similarly in the case of solids'' (presumably meaning for cubes and higher powers). No proof is given and no hints.  This implies that $\sqrt{2}$ is irrational and much more.   Such claims were evidently a surprise to the Pythagorean school (as in the probably apocryphal story of the demise of Hippasus).\footnote{Less apocryphally but a lot later, Pappus of Alexandria (290-350 AD) wrote that   ignorance  that incommensurables exist is ``a brutish and not a human state, and I am verily ashamed, not for myself only, but for ... what this whole generation believes, ... that commensurability is necessarily a quality of all magnitudes.''}
    So by the time of Euclid, irrational numbers were not just a theoretical and unlikely concept, but rather a practical one that mathematicians were studying, and entered Euclid's considerations in several different ways.

 In number theoretic terms,  Theaetetus claims that  if a positive integer is the $k$th power of a rational then it is the $k$th power of an integer. 
 This  follows  from results in Books 7, 8 and 10 of Euclid's \emph{Elements}:\footnote{However, Euclid struggles with his notation to explain these ideas so the proofs appear at times to meander, which  has even led Mazur (end of \S 6 of \cite{Maz}) to doubt that Euclid had a complete proof. I give his proof above, but gleaned from our earlier discussions from section 7 (extending 7.1 and 7.2 to  common measures of $k+1$ lengths),  \S  8.2 and its corollary, 8.3, 8.8  10.9 and 10.10 of the \emph{Elements}.}
   For Euclid an  integer $n$ is a $k$th power if there is a sequence of positive integers $a_0,\dots,a_k$ for which the ratios $a_{j+1}:a_j$ are all equal and $a_0:a_k=1:n$.
    Measuring $a_0$ and $a_1$ by their least common measure (that is, by dividing out by gcd$(a_0,a_1)$) we obtain $a_0:a_1=a:b$ where $a$ and $b$ are coprime positive integers.   The sequence $a^{k}, a^{k-1}b,\dots,b^k$ has the same ratios property as $a_0,\dots,a_k$ and the gcd of all of its elements is 1, so it can be obtained from 
 dividing $a_0,\dots,a_k$ by their least common measure.  In particular $a^k:b^k=1:n$ where gcd$(a^k,b^k)=$gcd$(1,n)=1$, and so $a^k=1$ and $b^k=n$.\footnote{Euclid evidently had no way to arrive at a $k$th power other than by multiplying together $k$ successive ratios (of lengths). If we go modern on this point then Euclid's proof is rather elegant: If $n^{1/k}$ is rational then we may write it as a reduced fraction $a/b$. Therefore $n/1=a^k/b^k$ with  $(n,1)=(a^k,b^k)=1$ so these are equal reduced fractions. However every positive rational is represented by just one reduced fraction (where the numerator and denominator are coprime positive integers) and so $a^k=1$ and $b^k=n$.}

Euclid's proof does not seem (to me) to be beyond what Theaetus might have done, though some historians are 
skeptical.\footnote{In the dialogue, Theaetus says: ``\emph{Theodorus here was drawing diagrams to show us something about powers ---  namely that the side lengths of a square of three square feet and also of one of five square feet aren't commensurable  with the side lengths of a square of one square foot; and so on, working out each case individually up to seventeen square feet at which
point, for some reason, he stopped}''. Even prior to Euclid, Greeks knew to work with squarefree areas (they knew that the area of a square is scaled by the square of the scaling of the length) and, since proving $\sqrt{2}$ is irrational was proved using an argument about parity in Aristotle's \emph{Prior Analytics}, historians guess that Theodorus therefore knew to avoid discussing squares of even square feet area by some more general parity argument.  Thus to make sense of the above we need an (unaccounted for) proof in which it is much tougher to prove $\sqrt{19}$ is irrational than $\sqrt{13}$; 
I have my doubts since Plato does not really exhibit any math skills nor shows understanding of the underlying ideas. Indeed Theaetus's main contribution in the dialogue seems to have been to invent adequate terminology to formulate the generalization to higher dimension, relative inanities which led to praise from Socrates (according to Plato). See also Hardy \& Wright, \cite{HW} sections 4.4 and 4.5 for their  attempt to reconstruct this history.}

   We have seen that if the lengths  of   two pieces of string are commensurable then the anthyphairesis process terminates after finitely many steps.\footnote{We actually proved that this is necessary as well as sufficient.}   Therefore if   the anthyphairesis process never terminates then the lengths of the original pieces of string are incommensurable,\footnote{This is proved by Euclid  in Book X, Proposition 2 where he uses the ``Eudoxus principle'' (that, if you keep on halving a quantity it gets arbitrarily small!) which could not have been a surprise to anyone (though Eudoxus had the genius to realize it needed proving), but the proof I just gave does not use this.} plausibly understood by mathematicians before Euclid.
Euclid's final touches on the anthyphairesis algorithm presents an impressively complete theory of commensurability.

  Hippasus is also credited with proving that the golden ratio $\phi:=\frac{\sqrt{5}+1}2$ is irrational, and I believe we can guess at his proof: Euclid defined the golden ratio as follows: Cut a length of string into two pieces. If the ratio of the length of the original piece of string to the largest piece, is the same as the ratio of the largest piece to the smallest piece, then that ratio is the golden ratio.\footnote{So if, after cutting, the two pieces of string have lengths $1>u$, say, then the original had length $1+u$, which we denote by $\phi$. Therefore 
  $1+u:1=1:u$ so that $u^2+u-1=0$  and $\phi^2-\phi-1=0$.} Now we apply anthyphairesis beginning with the original piece of the string and the larger part. The larger part can only be removed once from the original piece, and so we are left with the larger and smaller parts.
   Now, by definition,
    
 \centerline{Original Length : Largest Part  =  Largest Part : Smallest part}
 
 \noindent and so these two ratios have the same sequence of quotients (as in our previous discussion). But the second is obtained from the first by removing the quotient $1$ at the first step and, since the list of quotients are the same, therefore they must be 
 \[
 1,1,1,\cdots \text{ (ad infinitum).}
 \]
 This never terminates so the ratio, $\phi$, is incommensurable and therefore irrational.\footnote{This sequence of quotients shows that the continued fraction for $\phi$ is $[1,1,1,\dots]$.}

 With some effort one can create geometric proofs that $\sqrt{2}, \sqrt{3},\dots$ are irrational in the ancient Greek style, but each requires its own ad hoc argument. . For example, to prove that $\sqrt{2}$ is irrational we begin by putting four unit squares together (to create a square with side lengths 2, and area 4), then we draw the main diagonal of each the four squares to create an internal diamond square, of area 2 (since each diagonal halves the area of the four unit squares). The diagonals then give us lines of length $\sqrt{2}$, and cordon off a right-angled isoceles triangle. If 
 $2:1$ is commensurate  with a square like $a^2:b^2$ (where $a$ and $b$ are integers)  then we can scale up the side lengths, so the original squares have side lengths $a$, and the diagonal length $b$. By Pythagoras's theorem $b^2=2a^2$ which implies $b$ is even.  We now split the triangle in half, to again obtain a right-angled isoceles triangle but smaller, with hypoteneuse length $a$ and other sides of length $b/2$ (which is an integer since $b$ is even). There are now several ways to get to a contradiction and it is a matter of some debate which route the ancient Greeks took from here (perhaps by again using non-terminating anthyphairesis, or perhaps by noting that we could assume $(a,b)=1$ but we can deduce they must both be even, or ...).

Proving that $\pi$ is irrational remained elusive for more than two thousand years after the publication of the \emph{Elements}. It was only in the 
1760s that Lambert  proved that   $\pi$ is irrational, and then rather shockingly Lindemann proved that  $\pi$ is  transcendental in 1882.  There are many closely related questions that remain unresolved.

\section{The continuing saga of cultural bias}
The B\'ezout identity saga is to me a case of cultural misunderstanding. The modern era mainstream mathematical culture regards the   integers as the building blocks of mathematics, and for good reason; it is a very useful way to proceed. Euclid put us on that path by his genius in realizing that his results about measuring lengths hold independently of measuring. But today we scarcely recognize that as genius because we are educated to start with the integers and to go from there. 

Arguably how Euclid and his predecessors developed number theory, from measuring lengths, is more natural but less broadly applicable until you make the realizations that Euclid made.  We must not misinterpret his contributions based on our approach to thinking about numbers. This suggests we need to see that other ways to develop mathematics may have similarly validity and scope that is not obvious from our perspective.

We have persuaded ourselves that our perspective is so correct that it is commonly said and believed, without contradiction, that mathematics is a universal language, that 2+2 is always 4, so there can be no debate. However different cultures develop their tools, physical and mental, according to their needs and environment. Even breakthrough individuals, whether it be Euclid or Gauss, have tremendous insight based on understanding developments within their own time and culture.

Our culture has highly prioritized base 10 mathematics. But the  Inuit people in northern Quebec  learn numbers in a base 20 system,
probably because one can use 10 fingers and 10 toes to count,\footnote{This guess is  supported by the fact that the number 5 in Inuktitut derives from the word for ``hand'', 10 from ``top'' (for the upper body) and 20 from ``man is complete''.} and the transition between the two can be fraught with unexpected difficulties.
There have been many cultures, including ancient Mesopotamia and ancient Egypt, that have used base 12, probably 
because they counted with finger ``hinges'' on one hand, three on each of 4 fingers.\footnote{The internet pundits suggest that base 12 comes from the number of lunar cycles (months) in a year, but that seems like a far too sophisticated understanding for deciding upon the basis of a counting system.} This is reflected in how we measure time (seconds in a minute, minutes in an hour, and hours in a day), lengths (12 inches in a foot, $12^3$ UK yards in a mile), and eggs (by the ``dozen'').  The Telefol people in Papua New Guinea use base 27, involving many body parts in their basic counting, and who knows what else has worked in a given time and place.  It seems that all counting mechanisms are based on a  1-to-1 correspondence between the number of body parts of a certain type and a set of objects that needs to be measured.

For the Inuit there are several further challenges in participating in the broad discussion of mathematics since their language,
Inuit was an oral language until writing was introduced relatively recently. Therefore the numbers were not evidently separated from their uses:
 The digit ``3'' is \emph{pingasuk}; three objects is \emph{pingasut};
groups of three is  \emph{pingasuunaartitut};  a card numbered ``3'' is \emph{pingasulik};
three objects inside something else is \emph{pingasutalik}; and a pattern involving three objects is 
\emph{pingasuupaarqisimagusingha}. When we write these down we see that `pingasu'' is evidently the root, there in each usage, and with our background we immediately separate it in our mind to represent the concept of ``3''. But in the Inuit culture  how the 3 is used is important, and the number is integrated into a concept as it is used.  
Each of our digits (0 through to 10) has its own name but in Inuktitut, the digits began at ``3'' since there was no specific need to specify 1 or 2. Moreover the digits 6 to 9 are all composite expressions to explain the size, the most complex being
the digit ``7'', \emph{sitamaujunngigartut}, which translates as ``There are not exactly many fours'', ``many'' here meaning ``two''.
That can make it difficult for their children to quickly adapt to our more abstract approach to number, while our limited understanding about their use of number means we have had too little empathy for their challenges.  The subsequent introduction of the word ``\emph{atausik}'' for the digit ``1'' has added to the confusion -- it translates to ``indivisible'' which may then make it puzzling how one obtains fractions!

One geometric issue is that straight lines are a rarity in the  far north where traditional housing and other man made objects traditionally follow the contours of nature, so Inuits do not intuitively think of a ``line'' as straight. Indeed the translation of straight line is ``adopted line''.
To learn more about mathematics in this community, and indeed how many supposedly universal concepts are perceived differently, see \cite{Poi}. I should add that similar issues have been seen to arise in various cultures that have been historically isolated from citified society.

Nonetheless
Bishop \cite{bish} was able to identify six mathematical activities that seem to be found in every culture, though different cultures may develop these quite differently:

--- \emph{Counting}: Comparing and ordering sets of objects;

--- \emph{Spatial awareness}:  Explanations by way of drawings, models, symbols and words;

--- \emph{Measuring}: Using object and tools to quantify lengths, areas, volumes, and weights;

--- \emph{Design}:  Creation to help understanding, or for decoration;

--- \emph{Games}:  Particularly developing formal rules that need to be followed;

--- \emph{Exploration}:  Finding different ways to explain a phenomenon.

\noindent It is interesting to muse over how these might be developed differently from the developments that we know about.

\vskip 1.5in

\noindent For many of the French references in the bibliography the reader might consult  the freely available French national library, gallica.fr. 
For a  historian's perspective on the arithmetic in Euclid, see \cite{Sa}.

 \bibliographystyle{plain}

\begin{thebibliography}{99}



 \bibitem{Ba}  Claude-Gaspar Bachet, 
\emph{Probl\`emes plaisants et d\'electables  qui se font par les nombres},
Lyon 1612. (https://gallica.bnf.fr/ark:/12148/bpt6k8709271v)

 \bibitem{Ba2}  Claude-Gaspar Bachet, 
\emph{Probl\`emes plaisants et d\'electables  qui se font par les nombres}, 2nd ed.,
Lyon 1624. (https://gallica.bnf.fr/ark:/12148/bpt6k8702666r)

 \bibitem{bish} A.J. Bishop,
 \emph{Mathematics education in its cultural context}, in
Education studies in mathematics,
\textbf{19} (1988), 179-191. 
 

 \bibitem{Bo}  Nicolas Bourbaki, 
 \emph{\'El\'ements de math\'ematique: Alg\`ebre} 7 (1949).
 
 \bibitem{Bu}   M. F. Burnyeat,
  \emph{The philosophical sense of Theaetetus' mathematics},
 Isis,  69  (1978),  489--513.

 
 \bibitem{Co} Pierre Colmez, 
 \emph{\'El\'ements d'analyse et d`alg\`ebre (et de th\'eorie des nombres)},
 Les \'Editions de l'\'Ecole polytechnique (2011).
 
 
  \bibitem{Disq}  C.F. Gauss, 
\emph{Disquisitiones Arithmeticae}, (1804). (https://gallica.bnf.fr/ark:/12148/bpt6k3356j)

\bibitem{HW} G. H. Hardy  and  E. M. Wright,
\emph{An Introduction to the Theory of Numbers},
Oxford (1938).

\bibitem{Has} Brendon Hassett,
\emph{Introduction to algebraic geometry},
Cambridge (2012).

\bibitem{Haw}  Joel D. Hawkins, 
\emph{Proof and the Art of Mathematics},
MIT press (2023).

\bibitem{Hill} Richard Hill,
\emph{Introduction to number theory},
World scientific, (2019).


\bibitem{JJ} Gareth A. Jones, J. Mary Jones,
\emph{Elementary Number Theory},
Springer Undergraduate Mathematics Series, Springer London 1998.

\bibitem{Maz} Barry Mazur,
\emph{How did Theaetetus prove his theorem?} in 
``The Envisioned Life: Essays in Honor of Eva Brann'',
 Paul Dry Books (2007), pgs 227--250.
 
 
\bibitem{Poi} Louise Poirier,
\emph{Teaching mathematics and the Inuit community},
Canadian Journal of Science, Mathematics and Technology Education
\textbf{7} (2007), 53--67.

\bibitem{Sa}  Ken Saito, 
\emph{Re-examination of the Different Origins of the Arithmetical Books of Euclid's Elements},
Historia Mathematica 47 (2019), 39-53.

\bibitem{Sam} Pierre Samuel,
\emph{Th\'eorie alg\'ebrique des nombres},
Hermann, Paris, 1967.

\bibitem{Ten}   G\'erald Tenenbaum, 
\emph{Introduction to analytic and probabilistic number theory},
Grad. Stud. Math. \textbf{163} American Mathematical Society, Providence, RI, 2015. 

\bibitem{We}  Andr\'e Weil, 
\emph{Number Theory: An approach through history from Hammurapi to Legendre},
Birkh\"auser, New York (2001).


\end{thebibliography}

 \enddocument